\documentclass[3p,fleqn]{elsarticle}

\usepackage{amssymb}
\usepackage{amsbsy,amsfonts,amsmath,amssymb,bbm,calc,caption,color,dsfont,graphics}
\usepackage{graphicx,ifthen,latexsym,mathrsfs,multicol}
\usepackage{natbib,overpic,pifont,psfrag,rotating,stmaryrd,theorem}

\newtheorem{theorem}{Theorem}
\newtheorem{lemma}[theorem]{Lemma}

\newdefinition{remark}{Remark}
\newdefinition{corollary}{Corollary}

\numberwithin{equation}{section}

\numberwithin{theorem}{section}
\journal{}

\begin{document}
\begin{frontmatter}
\title{General constraint preconditioning iteration method for singular saddle-point problems}

\author[lzu]{Ai-Li Yang\corref{yal}}\ead{yangaili@lzu.edu.cn; cmalyang@gmail.com}
\cortext[yal]{Corresponding author. Tel.: +86 931 8912483; fax: +86 931 8912481.}
\author[lzu]{Guo-Feng Zhang}\author[lzu,fup]{Yu-Jiang Wu}

\address[lzu]{School of Mathematics and Statistics, Lanzhou University, Lanzhou 730000, PR China}
\address[fup]{Department of Mathematics, Federal University of Paran\'{a}, Centro Polit\'{e}cnico, CP: 19.081, 81531-980, Curitiba, PR
Brazil}

\begin{abstract}
For the singular saddle-point problems with nonsymmetric positive definite $(1,1)$ block, we present a general constraint preconditioning (GCP) iteration method based on a singular constraint preconditioner. Using the properties of the Moore-Penrose inverse, the convergence properties of the GCP iteration method are studied. In particular, for each of the two different choices of the $(1,1)$ block of the singular constraint preconditioner, a detailed convergence condition is derived by analyzing the spectrum of the iteration matrix. Numerical experiments are used to illustrate the theoretical results and examine the effectiveness of the GCP iteration method. Moreover, the preconditioning effects of the singular constraint preconditioner for restarted generalized minimum residual (GMRES) and quasi-minimal residual (QMR) methods are also tested.

\end{abstract}
\end{frontmatter}

\section{Introduction}\label{sec:1}
Consider the following large, sparse singular saddle-point problems
\begin{equation}\label{01}
A\,x:=\left( \begin{array}{cc} W & B^T \\  -B & 0 \end{array} \right)
\left( \begin{array}{c} u \\ v \end{array} \right)
=\left( \begin{array}{c} f \\ g \end{array} \right)=b,
\end{equation}
where $W\in\mathbb{R}^{n\times n}$ is nonsymmetric positive definite and $B\in\mathbb{R}^{m\times n}$ is rank deficient, i.e., $\text{rank}(B)<m\leq n$, $b\in \mathbb{R}^{n+m}$ is a given vector in the range of saddle-point matrix $A\in\mathbb{R}^{(n+m)\times(n+m)}$. Such kind of linear systems arise in many application areas, such as computational fluid dynamics, computational genetics,
mixed finite element approximation of elliptic partial differential equations, constrained optimization, optimal control, weighted least-squares problems, electronic networks, computer graphics etc; see \cite{BPW20051,SSY19981,WSY2004139,ZBY2009808} and references therein.

When the saddle-point matrix $A$ in \eqref{01} is nonsingular, which requires $B$ being of
full row rank, a number of iteration methods and preconditioning techniques have been proposed to approximate the unique solution of the nonsingular saddle-point problem \eqref{01}; see \cite{BGP20041,BGL20051,EG19941645,FRSW1998527,KGW20001300,B2009447,BW20082900} and references therein.
Within these results, the constraint preconditioner of the form
\begin{equation}\label{02}
M=\left( \begin{array}{cc} P & B^T \\  -B & 0 \end{array} \right),
\end{equation}
with $P$ being positive definite was widely analyzed; see \cite{KGW20001300,BG20051125,Cao2002121,GW1998530,BNW2009410}. Based on this preconditioner $M$, Golub and Wathen \cite{GW1998530} studied the following basic iteration scheme
\begin{equation}\label{32}
    x^{(k+1)}=x^{(k)}+M^{-1}(b-Ax^{(k)}).
\end{equation}
We call this scheme constraint preconditioning iteration method if $M$ is chosen to be the nonsingular constraint preconditioner \eqref{02}.
Let $H$ and $S$ be respectively the symmetric and the skew-symmetric parts of matrix $W$, i.e.,
\[W=H+S,\quad\text{with }~ H=\frac12(W+W^T)~\text{ and }~S=\frac12(W-W^T).\]
The convergence properties of \eqref{32} were derived by Golub and Wathen \cite{GW1998530} when matrix $P$ in $M$ is chosen to be a multiple of the symmetric part of $W$, i.e., $P=\omega H$ with $\omega>0$. If $W$ is not far from a symmetric matrix (i.e., $\|S\|/\|H\|$ is a small number), the preconditioner $M$ with $P=\omega H$ is very efficient. However, as can be expected, performance of this preconditioner with symmetric $P$ deteriorates when $W$ is essentially nonsymmetric ($\|S\|/\|H\|\approx 1$ or larger). To overcome this deficiency, Botchev and Golub \cite{BG20051125} proposed a novel constraint preconditioner by choosing the $(1,1)$ block of \eqref{02} with
\begin{equation}\label{33}
  P=\frac1{\omega}(I+\omega L_s)(I+\omega U_s),
\end{equation}
where $\omega>0$, $L_s$ and $U_s$ are, respectively, lower and upper triangular parts of the matrix $S$ satisfying $L_s+U_s=S$ and $U_s=-L_s^T$.
The preconditioner $M$ with the new choice of $P$ used for the iteration scheme \eqref{32} was proved to be efficient and robust for solving nonsingular saddle-point problems \eqref{01} with $W$ being nonsymmetric. Moreover, as a preconditioner, it also can improve the convergence rate of GMRES method.


When $B$ is rank deficient, both of the saddle-point matrix $A$ in \eqref{01} and the constraint preconditioner $M$ in \eqref{02} are singular. The linear systems \eqref{01} are called as singular saddle-point problems. Some authors have studied iterative methods or preconditioners for this kind of singular problems and obtained many important and interesting results; see \cite{WSY2004139,ZBY2009808,Bai2010171,BaiWangYuan200375,ZhangShen2013116,CZ2013doi,ZLW2014334}. Owing to the singularity of matrix $M$, iteration scheme \eqref{32} can not be used to solve singular saddle-point problems \eqref{01}. In 2008, Cao \cite{Cao20081382} proposed an iteration scheme by replacing $M^{-1}$ with $M^\dag$ in \eqref{32} to solve general singular linear systems $Ax=b$, that is
\begin{equation}\label{03}
  x^{(k+1)}=x^{(k)}+M^\dag (b-Ax^{(k)}),
\end{equation}
where $M$ is a singular matrix depending on the coefficient matrix $A$, $M^\dag$ is the Moore-Penrose inverse of matrix $M$ satisfying the following Moore-Penrose equations:
\begin{equation}\label{11}
  MM^\dag M=M,\quad (M^\dag M)^*=M^\dag M,\quad (MM^\dag)^*=MM^\dag,\quad M^\dag MM^\dag=M^\dag.
\end{equation}
Iteration scheme \eqref{03} was used later to solve the range-Hermitian singular linear systems by Zhang and Wei in \cite{ZhangWei2010139}, the numerical efficiencies of this method were also verified. We call iteration scheme \eqref{03} the general constraint preconditioning (GCP) iteration method if $M$ is a singular constraint preconditioner of the form \eqref{02}.

In this work, we are especially interested in the case that matrix $B$ is rank deficient, which means the saddle-point matrix $A$ in \eqref{01} and the constraint preconditioner $M$ in \eqref{02} are both singular. We use GCP iteration method \eqref{03} to solve the singular saddle-point problems \eqref{01}. The remainder part of this work is organized as follows. In Section \ref{sec:2}, we give the convergence properties of GCP iteration method \eqref{03} with $M$ being of the form \eqref{02} and $P$ being any positive definite matrix. For each of the two different choices of the matrix $P$, i.e., $P=\omega H$ and $P=(1/\omega)(I+\omega L_s)(I+\omega U_s)$, a detailed condition that guarantees the convergence of the GCP iteration method is derived in Section \ref{sec:3}. In Section \ref{sec:4}, numerical results show that the GCP iteration method \eqref{03}, no matter as a solver or as a preconditioner for GMRES(10) and QMR methods, is robust and efficient.
Finally in Section \ref{sec:5}, we end this work with a brief conclusion.


\section{Convergence properties}\label{sec:2}
In this section, we analyze the convergence properties of the GCP iteration method \eqref{03} with $M$ being defined in \eqref{02} and $P$ being positive definite (maybe not symmetric). First, we present the following convergence result of iteration scheme \eqref{03} with any singular matrix $M$:
\begin{lemma}\label{04}\cite{Cao20081382}
Iteration scheme \eqref{03} is convergent if and only if the following three conditions are fulfilled:
  \begin{enumerate}
    \item $\text{null}\,(M^\dag A)=\text{null}\,(A)$;
    \item index$\,(I-T)=1$, or equivalently, rank$\,(I-T)=$rank$\,((I-T)^2)$, where $T:=I-M^\dag A$ is the iteration matrix of \eqref{03};
    \item $\gamma(T)=\max\{|\lambda|:\lambda\in\sigma(T)\backslash\{1\}\}<1$, where $\sigma(T)$ is the spectral set of matrix $T$.
  \end{enumerate}
\end{lemma}

In the following subsections, we analyze the convergence properties of GCP iteration method \eqref{03}, i.e., $M$ is singular matrix of the form \eqref{02}, according to the three conditions of Lemma \ref{04}.
\subsection{The first condition of Lemma \ref{04}}
For any $x\in \mathbb{R}^{n+m}$ satisfying $MM^\dag Ax=0$, we have $M^\dag Ax=M^\dag (MM^\dag Ax)=0$, which gives $\text{null}(MM^\dag A)\subseteq\text{null}(M^\dag A)$. Since $\text{null}(MM^\dag A)\supseteq\text{null}(M^\dag A)$ is obvious, we can obtain that
\begin{equation}\label{38}
\text{null}(MM^\dag A)=\text{null}(M^\dag A).
\end{equation}
From the definition of Moore-Penrose inverse, matrix $M^\dag$ can be written as \cite{ZhangShen2013116}:
\begin{equation}\label{05}
    M^\dag=\left(
             \begin{array}{cc}
               P^{-1}-P^{-1}B^TE^\dag BP^{-1} & -P^{-1}B^TE^\dag \\
               E^\dag BP^{-1} & E^\dag\\
             \end{array}
           \right),
\end{equation}
where $E=BP^{-1}B^T$. Owing to $EE^\dag B=E^\dag EB=B$ \cite{ZhangShen2013116}, it follows that
\[MM^\dag=\left( \begin{array}{cc}  I & 0 \\  0 & EE^\dag  \end{array} \right)\]
and
\begin{equation}\label{07}
MM^\dag A = \left( \begin{array}{cc} I & 0 \\  0 & EE^\dag \end{array} \right)
            \left( \begin{array}{cc} W & B^T \\ -B & 0 \end{array} \right)
          = \left( \begin{array}{cc} W & B^T \\ -EE^\dag B & 0 \end{array} \right)
          = A.
\end{equation}
Thus, using \eqref{38} and \eqref{07}, we finally obtain that $\text{null}(M^\dag A)=\text{null}(MM^\dag A)=\text{null}(A)$.

\subsection{The second condition of Lemma \ref{04}}
Since $T=I-M^\dag A$, the second condition of Lemma \ref{04} holds if $\text{null}((M^\dag A)^2)=\text{null}(M^\dag A)$. Owing to $\text{null}((M^\dag A)^2)\supseteq\text{null}(M^\dag A)$, we only need to prove $\text{null}((M^\dag A)^2)\subseteq\text{null}(M^\dag A)$ in the following.

From \eqref{05}, matrix $M^\dag A$ can be written as
\begin{equation}\label{06}
M^\dag A=\left( \begin{array}{cc}  P^{-1}W-P^{-1}B^TE^\dag BP^{-1}W+P^{-1}B^TE^\dag B & 0 \\  E^\dag BP^{-1}W-E^\dag B & E^\dag E \end{array} \right).
\end{equation}
Let $x=\left( x_1^T, x_2^T\right)^T\in \mathbb{R}^{n+m}$, with $x_1\in\mathbb{R}^{n}$ and $x_2\in\mathbb{R}^{m}$, satisfying $(M^\dag A)^2x=0$. Denote
\begin{equation}\label{08}
    M^\dag Ax=\left(
                  \begin{array}{c}
                    \left(P^{-1}W-P^{-1}B^TE^\dag BP^{-1}W+P^{-1}B^TE^\dag B\right)x_1 \\
                    E^\dag BP^{-1}Wx_1-E^\dag Bx_1+ E^\dag Ex_2
                  \end{array}
                \right)=:\left(\begin{array}{c} y_1 \\  y_2 \end{array} \right)=y.
\end{equation}
In the following, we prove $M^\dag Ax=y=0$. Using $M^\dag Ay=(M^\dag A)^2x=0$ and $\text{null}(M^\dag A)=\text{null}(A)$ gives $Ay=0$, i.e.,
\begin{equation}\label{09}
Wy_1+B^Ty_2=0\quad\text{and}\quad -By_1=0.
\end{equation}
Since $W$ is nonsingular, solving $y_1$ from the first equality of \eqref{09} and taking it into the second equality, we have $BW^{-1}B^Ty_2=0$.
Hence, it follows that
\[(B^Ty_2)^TW^{-1}(B^Ty_2)=y_2^T(BW^{-1}B^Ty_2)=0.\]
Owing to the positive definiteness of matrix $W^{-1}$, we further obtain that $B^Ty_2=0$. Taking it into the first equality of \eqref{09} and noticing that $W$ is nonsingular, we first get $y_1=0$.

From \eqref{08}, we have with $E=BP^{-1}B^T$ and $EE^\dag B=B$ that
\[
0=By_1=BP^{-1}Wx_1-BP^{-1}B^TE^\dag BP^{-1}Wx_1+BP^{-1}B^TE^\dag B x_1=Bx_1.
\]
Thus, vector $y_2$ in \eqref{08} can be written as
\[y_2=E^\dag BP^{-1}Wx_1+E^\dag E x_2.\]
Using $B^Ty_2=0$ and $B^TE^\dag E=(E^\dag EB)^T=B^T$ gives
\[
\begin{split}
  0=E^\dag BP^{-1}(B^Ty_2)=&E^\dag BP^{-1}(B^TE^\dag BP^{-1}Wx_1+B^TE^\dag E x_2)\\
  =&E^\dag BP^{-1}(B^TE^\dag BP^{-1}Wx_1+B^T x_2)\\
  =&E^\dag BP^{-1}Wx_1+E^\dag E x_2=y_2.
\end{split}
\]
Therefore, we obtain $y=(y_1^T,y_2^T)^T=0$. The GCP iteration method \eqref{03} satisfies the second condition of Lemma \ref{04}.

\subsection{The third condition of Lemma \ref{04}}
From \eqref{06}, we have
\begin{equation}\label{10}
I-M^\dag A=\left( \begin{array}{cc}  X(P-W) & 0 \\  -E^\dag BP^{-1}W+E^\dag B & I-E^\dag E \end{array} \right),
\end{equation}
where $X=P^{-1}-P^{-1}B^TE^\dag BP^{-1}$. Let the singular value decomposition of $E$ be
\[E=U\left( \begin{array}{cc}\Sigma & 0 \\  0 & 0 \\ \end{array}\right)V^T,\]
where $U$ and $V$ are two orthogonal matrices. Then, the $(2,2)$ block of matrix $I-M^\dag A$ in \eqref{10} becomes
\[
I-E^\dag E=I-V\left( \begin{array}{cc} \Sigma^{-1} & 0 \\ 0 & 0 \\  \end{array} \right)U^TU
\left( \begin{array}{cc} \Sigma & 0 \\ 0 & 0 \end{array}\right)V^T
=V\left(\begin{array}{cc} 0 & 0 \\ 0 & I  \end{array} \right)V^T.
\]
Hence, matrix $I-E^\dag E$ only has two different eigenvalues, which are $0$ and $1$. Therefore, we have
\[\gamma(I-M^\dag A)=\gamma(X(P-W)).\]

Using Lemma \ref{04}, we summarize this section with the following convergence result of GCP iteration method \eqref{03}.
\begin{theorem}\label{19}
  For the singular saddle-point problems \eqref{01}, let matrices $W$ and $B$ be positive definite and rank deficient, respectively. Then, the GCP iteration scheme \eqref{03}, with $M$ defined in \eqref{02} and $P$ being positive definite, is convergent if and only if
  \[\gamma(X(P-W))<1,\]
  where $X=P^{-1}-P^{-1}B^TE^\dag BP^{-1}$.
\end{theorem}

\section{The choices of matrix $P$}\label{sec:3}
In this section, based on the two different choices of submatrix $P$ of singular constraint preconditioner $M$, we further analyze the convergence properties of the GCP iteration methods \eqref{03}.
\subsection{Symmetric submatrix $P$}
In this subsection, we suppose that $P$ is symmetric positive definite. Thus, $X(P-W)$ is similar to
\begin{equation}\label{16}
P^{1/2}X(P-W)P^{-1/2}=P^{1/2}XP^{1/2}\left(I-P^{-1/2}WP^{-1/2}\right).
\end{equation}
For the matrix $P^{1/2}XP^{1/2}$, we have the following properties:
\begin{lemma}\label{12}
  Let $P$ be symmetric positive definite. Then, the $n\times n$ real matrix $P^{1/2}XP^{1/2}$ is symmetric, which has $n-\text{rank}(B)$ eigenvalues $\lambda=1$ and $\text{rank}(B)$ eigenvalues $\lambda=0$.
\end{lemma}
\emph{Proof}. The symmetry of $P^{1/2}XP^{1/2}$ is obvious, we only need to prove the remainder part of this lemma. Using the definition of matrix $X$ gives
\begin{equation}\label{13}
  P^{1/2}XP^{1/2}=I-P^{-1/2}B^TE^\dag BP^{-1/2}.
\end{equation}
Since $E=BP^{-1}B^T$ and $E^\dag EE^\dag=E^\dag$, it follows
\[(P^{-1/2}B^TE^\dag BP^{-1/2})^2=P^{-1/2}B^TE^\dag BP^{-1/2},\]
which means $P^{-1/2}B^TE^\dag BP^{-1/2}$ is a projection matrix.
Hence, from \eqref{13}, matrix $P^{1/2}XP^{1/2}$ is also a projection. The eigenvalues of symmetric matrix $P^{1/2}XP^{1/2}$ are $0$ or $1$.

Denoting $G=BP^{-1/2}$, matrix $P^{-1/2}B^TE^\dag BP^{-1/2}$ can be rewritten as
\begin{equation}\label{39}
P^{-1/2}B^TE^\dag BP^{-1/2}=G^T(GG^T)^\dag G.
\end{equation}
Since $\text{rank}(G)=\text{rank}(B)$, it follows that
\begin{equation}\label{14}
  \text{rank}(G^T(GG^T)^\dag G)\leq\text{rank}(G)=\text{rank}(B)
\end{equation}
and
\begin{equation}\label{15}
  \text{rank}(G^T(GG^T)^\dag G)\geq\text{rank}(GG^T(GG^T)^\dag GG^T)=\text{rank}(GG^T)=\text{rank}(G)=\text{rank}(B).
\end{equation}
Combining \eqref{14} and \eqref{15} and using \eqref{39}, we have
\[
\text{rank}(P^{-1/2}B^TE^\dag BP^{-1/2})=\text{rank}(G^T(GG^T)^\dag G)=\text{rank}(B).
\]
Inasmuch as $P^{-1/2}B^TE^\dag BP^{-1/2}$ is a projection matrix, it has $\text{rank}(B)$ eigenvalues $\lambda=1$ and $n-\text{rank}(B)$ eigenvalues $\lambda=0$.
Using \eqref{13}, we finally obtain that matrix $P^{1/2}XP^{1/2}$ has $n-\text{rank}(B)$ eigenvalues $\lambda=1$ and $\text{rank}(B)$ eigenvalues $\lambda=0$. $~~~~~~~\Box$

When the symmetric part of $W$ is dominant, a particular choice of $P$ is $P=\omega H$, where $\omega$ is a real and positive constant and $H$ is the symmetric part of $W$. From Lemma \ref{12}, there exists an orthogonal matrix $Q$, such that
\[
P^{1/2}XP^{1/2}=Q\left( \begin{array}{cc}I_r & 0 \\  0 & 0 \\ \end{array}\right)Q^T,\quad\text{with } r=n-\text{rank}(B).
\]
Hence, the eigenvalues of $P^{1/2}XP^{1/2}\left(I-P^{-1/2}WP^{-1/2}\right)$ in \eqref{16} are the eigenvalues of
\begin{equation}\label{17}
    \left( \begin{array}{cc}I_r & 0 \\  0 & 0 \\ \end{array}\right)Q^T\left(I-P^{-1/2}WP^{-1/2}\right)Q.
\end{equation}
In other words, the nonzero eigenvalues of the matrix $P^{1/2}XP^{1/2}\left(I-P^{-1/2}WP^{-1/2}\right)$ are the eigenvalues of the $r\times r$ leading principle submatrix of
\begin{equation}\label{18}
    Q^T\left(I-P^{-1/2}WP^{-1/2}\right)Q=\left(1-\frac 1{\omega}\right)I-\frac 1{\omega}Q^TH^{-1/2}SH^{-1/2}Q.
\end{equation}
Let $i$ denote the imaginary unit, then matrix $iQ^TH^{-1/2}SH^{-1/2}Q$ is Hermitian. From the interlace theorem \cite{GolubVL1996,Wilkinson1988}, the eigenvalues $i\eta$ of the $r\times r$ leading principle submatrix of $Q^TH^{-1/2}SH^{-1/2}Q$ satisfy
\[
|\eta|\leq\rho(Q^TH^{-1/2}SH^{-1/2}Q)=\rho(H^{-1/2}SH^{-1/2}),
\]
where $\rho(\cdot)$ denotes the spectral radius of a matrix.
Since $Q^TH^{-1/2}SH^{-1/2}Q$ is skew-symmetric, the $r\times r$ leading principle submatrix of \eqref{18} is of the form
$\left(1-1/{\omega}\right)I-(1/{\omega})K$ with $K=-K^T$.
Hence, the eigenvalues of the $r\times r$ leading principle submatrix of \eqref{18} are of the form
$1-1/{\omega}-i\eta/\omega$ satisfying $|\eta|\leq\rho(H^{-1/2}SH^{-1/2})$.

Therefore, the eigenvalues of matrix $P^{1/2}XP^{1/2}(I-P^{-1/2}WP^{-1/2})$ (or $X(P-W)$) with $P=\omega H$ are either zero or else are of the form $1-1/{\omega}-i\eta/\omega$ with $\eta\in\mathbb{R}$ and $|\eta|\leq\rho(H^{-1/2}SH^{-1/2})$. We finally obtain that
\[\gamma(X(P-W))\leq \frac1{\omega}\left((\omega-1)^2+\rho\left(H^{-1/2}SH^{-1/2}\right)^2\right)^{1/2}.\]
Now, using Theorem \ref{19}, the following convergence result is obtained.
\begin{theorem}
  For the singular saddle-point problems \eqref{01}, let matrices $W$ and $B$ be positive definite and rank deficient, respectively. Then, the GCP iteration scheme \eqref{03}, with $M$ defined in \eqref{02} and $P=\omega H$, is convergent if
  \[\omega>\frac12\left(1+\rho\left(H^{-1/2}SH^{-1/2}\right)^2\right).\]
\end{theorem}

\subsection{Non-symmetric submatrix $P$}
In this subsection, we analyze the convergence properties of the GCP iteration method \eqref{03} with $M$ defined in \eqref{02} and $P$ defined in \eqref{33}.

First, we give the following property of matrix $P$:
\begin{lemma}\label{31}\cite{BG20051125}
  The matrix $P$ in \eqref{33} is positive definite if and only if $\omega<1/{\|L_s\|_2}$.
\end{lemma}

In the remainder of this subsection, we suppose that $\omega<1/{\|L_s\|_2}$, i.e., matrix $P$ in \eqref{33} is positive definite. From Theorem \ref{19}, the analysis of convergence properties of iteration scheme \eqref{03} is reduced to estimate the pseudospectral radius $\gamma(X(P-W))$. Let $P_H$ be the symmetric part of matrix $P$ in \eqref{33}, we can bound $\gamma(X(P-W))$ as
\begin{equation}\label{22}
    \gamma(X(P-W))=\gamma\left(P_H^{1/2}X(P-W)P_H^{-1/2}\right)\leq\left\|P_H^{1/2}XP_H^{1/2}\right\|_2\cdot\left\|P_H^{-1/2}(P-W)P_H^{-1/2}\right\|_2.
\end{equation}
Denote
\[
X_\delta=P^{-1}-P^{-1}B^TE^T(EE^T+\delta^2I)^{-1}BP^{-1}.
\]
It is easy to verify that $X=\lim_{\delta\rightarrow 0} X_\delta$, since $E^\dag=\lim_{\delta\rightarrow 0}E^T(EE^T+\delta^2I)^{-1}$ \cite{Albert1972}. From Sherman-Morrison-Woodbury formula \cite{GolubVL1996}, matrix $X_\delta$ can be rewritten as
\begin{equation}\label{23}
    X_\delta=\left(P+\frac1{\delta^2}B^TE^TB\right)^{-1}.
\end{equation}
Let $P_S$ be the skew-symmetric part of matrix $P$ in \eqref{33}, it follows that
\[
P_H^{1/2}X_\delta P_H^{1/2}=\left(I+\widetilde{P}_S+\frac1{\delta^2}L^TP^{-T}L\right)^{-1},
\]
where $\widetilde{P}_S=P_H^{-1/2}P_SP_H^{-1/2}$ is a skew-symmetric matrix and $L=B^TBP_H^{-1/2}$. Simple calculation gives
\begin{equation}\label{24}
\begin{split}
  \left\|P_H^{1/2}X_\delta P_H^{1/2}\right\|_2^2&=\max_{\|x\|_2=1}\frac1{\left((I+\widetilde{P}_S+\frac1{\delta^2}L^TP^{-T}L)x,(I+\widetilde{P}_S+\frac1{\delta^2}L^TP^{-T}L)x\right)}\\
  &\leq\frac1{1+\displaystyle\frac1{\delta^2}\min_{\|x\|_2=1}\left(L^T(P^{-T}+P^{-1})Lx,x\right)+\left\|\left(\widetilde{P}_S+\frac1{\delta^2}L^TP^{-T}L\right)x\right\|_2^2}.
\end{split}
\end{equation}
Since $P$ is positive definite, $\left(L^T(P^{-T}+P^{-1})Lx,x\right)\geq0$, which gives from \eqref{24} that
$\|P_H^{1/2}X_\delta P_H^{1/2}\|_2\leq 1$.
Let $\delta\rightarrow 0$, we have
\begin{equation}\label{25}
\left\|P_H^{1/2}X P_H^{1/2}\right\|_2=\lim_{\delta\rightarrow 0}\left\|P_H^{1/2}X_\delta P_H^{1/2}\right\|_2\leq 1.
\end{equation}

In the following, we estimate $\|P_H^{-1/2}(P-W)P_H^{-1/2}\|_2$. Inasmuch as $P-W=P_H-H$ is symmetric, it follows
\begin{equation}\label{28}
    \left\|P_H^{-1/2}(P-W)P_H^{-1/2}\right\|_2=\left\|I-P_H^{-1/2}HP_H^{-1/2}\right\|_2=\rho\left(I-P_H^{-1/2}HP_H^{-1/2}\right).
\end{equation}
Since matrix $H$ is symmetric positive definite, the eigenvalues of $I-P_H^{-1/2}HP_H^{-1/2}$ are inside the interval $(-1,1)$ if and only if
\[2(P_Hx,x)>(Hx,x),\quad\forall ~x\in\mathbb{R}^n\text{ and }x\neq 0,\]
or, say
$2P_H-H=(2/{\omega})I-2\,\omega L_sL_s^T-H$
is positive definite. Let $\lambda_{\max}(H)$ be the maximum eigenvalue of matrix $H$. Note that $\lambda_{\max}(L_sL_s^T)=\|L_s\|_2^2$, matrix $2P_H-H$ is positive definite if
\begin{equation}\label{26}
    \frac2{\omega}-2\,\omega\|L_s\|_2^2-\lambda_{\max}(H)>0.
\end{equation}
Since $\omega>0$, inequality \eqref{26} holds if and only if
\begin{equation}\label{27}
    0<\omega<\frac{-\lambda_{\max}(H)+\sqrt{\lambda_{\max}(H)^2+16\|L_s\|_2^2}}{4\|L_s\|_2^2}.
\end{equation}
Therefore, under condition \eqref{27}, we obtain from \eqref{28} that
\begin{equation}\label{29}
\left\|P_H^{-1/2}(P-W)P_H^{-1/2}\right\|_2<1.
\end{equation}
Combining \eqref{22}, \eqref{25} and \eqref{29}, we know that $\gamma(X(P-W))<1$ if $\omega$ satisfies condition \eqref{27}.
Recall Theorem \ref{19}, the convergence of GCP iteration scheme \eqref{03} can be established if we have another condition, i.e., matrix $P$ is positive definite.
Therefore, using Lemma \ref{31} and noticing that
\[\frac{-\lambda_{\max}(H)+\sqrt{\lambda_{\max}(H)^2+16\|L_s\|_2^2}}{4\|L_s\|_2^2}\leq\frac1{\|L_s\|_2},\]
we finally derive the following convergence results of GCP iteration scheme \eqref{03}:
\begin{theorem}
  For the singular saddle-point problems \eqref{01}, let matrices $W$ and $B$ be positive definite and rank deficient, respectively. Then, the GCP iteration scheme \eqref{03}, with $M$ defined in \eqref{02} and $P=(1/{\omega})(I+\omega L_s)(I+\omega U_s)$, is convergent if
  \[0<\omega<\frac{-\lambda_{\max}(H)+\sqrt{\lambda_{\max}(H)^2+16\|L_s\|_2^2}}{4\|L_s\|_2^2}.\]
\end{theorem}

\section{Numerical experiments}\label{sec:4}
In this section, we assess the feasibility and robustness of the GCP iteration method \eqref{03} with matrix $M$ being defined in \eqref{02} and $P$ being positive definite. In addition, the preconditioning effects of the singular constraint preconditioners $M$ for GMRES(10) and QMR methods will also be tested.

%

Consider the linearized version of the steady-state Navier-Stokes equations, i.e., the Oseen equations of the following form
\begin{equation}\label{34}
    \left\{
    \begin{split}
      &-\nu \Delta \textbf{u}+(\textbf{w}\cdot \nabla)\textbf{u}+\nabla p = \textbf{f},\quad \text{in }\Omega,\\
      &-\nabla\cdot \textbf{u}=0,\quad \text{in }\Omega,
    \end{split}\right.
\end{equation}
where $\Omega$ is an open bounded domain in $\mathbb{R}^2$, vector $\textbf{u}$ represents the velocity in $\Omega$, function $p$ represents pressure, and
the scalar $\nu>0$ is the viscosity constant. The test problem is a leaky-lid driven cavity problem in square domain $\Omega=(0,1)\times(0,1)$ with the wind field $\textbf{w} = (a(x,y), b(x,y))^T$ being
chosen as $a(x,y)= 8x(x-1)(1-2y)$ and $b(x,y)= 8y(2x-1)(y-1)$. The boundary conditions are $\textbf{u} = (0,0)^T$ on the three
fixed walls $(x = 0, y = 0, x = 1)$, and $\textbf{u} = (1, 0)^T$ on the moving wall $(y = 1)$.

Dividing $\Omega$ into a uniform $l\times l$ grid with mesh size $h=1/l$ and discretizing \eqref{34} by the "marker and cell" (MAC) finite difference scheme \cite{HW19652182,Elman19991299}, the singular saddle-point system \eqref{01} is obtained, where
\[
W=\left( \begin{array}{cc} F_1 & 0 \\ 0 & F_2 \end{array} \right)\in\mathbb{R}^{2l(l-1)\times 2l(l-1)},\quad B=(B_1,B_2)\in\mathbb{R}^{l^2\times 2l(l-1)},
\]
and
\[
F_i=\nu A_i+N_i\in\mathbb{R}^{l(l-1)\times l(l-1)},\quad i=1,2.
\]
The coefficient matrix $A$ of \eqref{01} has the following properties: $W$ is nonsymmetric and positive definite, $\text{rank}(B) = l^2-1$, thus $A$ is singular.

\begin{table}[h]
\caption{Choices of the preconditioner $M$} \label{tab:01}
\begin{center}
\begin{tabular}{*{2}{l}}
\hline\noalign{\smallskip}
Case no.& Preconditioner $M$ \\
\hline\noalign{\smallskip}
I  & $M$ defined in \eqref{02} with $P=\omega H$ \\
II  & $M$ defined in \eqref{02} with $P=(1/\omega)(I+\omega L_s)(I+\omega U_s)$ \\
III  & $M_b$ defined in \eqref{35} with $P=\omega H$\\
IV  & $M_b$ defined in \eqref{35} with $P=(1/\omega)(I+\omega L_s)(I+\omega U_s)$\\
V  &  $M_t$ defined in \eqref{36} with $P=\omega H$\\
VI  &  $M_t$ defined in \eqref{36} with $P=(1/\omega)(I+\omega L_s)(I+\omega U_s)$\\
\noalign{\smallskip}\hline
\end{tabular}
\end{center}
\end{table}

The efficiency of iteration scheme \eqref{03} will be tested by comparing its iteration
steps (denoted as IT), elapsed CPU time in seconds (denoted as CPU) with those of the iteration scheme \eqref{32}. For iteration scheme \eqref{03}, besides the singular constraint preconditioners $M$ defined in \eqref{02} with $P=\omega H$ and $P=(1/{\omega})(I+\omega L_s)(I+\omega U_s)$, we choose another kind of singular block diagonal preconditioners as
\begin{equation}\label{35}
    M_b=\left(
      \begin{array}{cc}
        P & 0 \\
        0 & BP^{-1}B^T \\
      \end{array}
    \right).
\end{equation}
For iteration scheme \eqref{32}, the following nonsingular block triangular preconditioners will be tested \cite{Elman19991299}:
\begin{equation}\label{36}
    M_t=\left(
      \begin{array}{cc}
        P & B^T \\
        0 & \frac1{\nu}h^2I \\
      \end{array}
    \right).
\end{equation}
The matrix $P$ in each of the preconditioners \eqref{35} and \eqref{36} is also chosen to be $P=\omega H$ and $P=(1/{\omega})(I+\omega L_s)(I+\omega U_s)$, respectively. The detailed test cases can be seen in Table \ref{tab:01}.

In the implementations, the iteration methods are started from zero vector and terminated once the current iterate $x^{(k)}$ satisfies
\begin{equation}\label{40}
\text{RES}=\frac{\|b-Ax^{(k)}\|_2}{\|b\|_2}< 10^{-6}.
\end{equation}
In addition, all codes were run in MATLAB [version 7.10.0.499 (R2010a)] in double precision and all experiments were performed on a personal computer with 3.10GHz central processing unit [Intel(R) Core(TM) Duo i5-2400] and 3.16G memory.

For the parameters $\omega$ in matrices $P$, we choose the experimentally found optimal ones, which result in the least iteration steps for iteration schemes \eqref{32}, \eqref{03} and preconditioned GMRES(10) and QMR methods; see Tables \ref{tab:02}, \ref{tab:03} and \ref{tab:04}. The sign "$-$" in the three tables is used to denote that the methods do not converge within IT$_{\max}=5000$ iteration steps for any $\omega\in(0,2000]$. Numerical experiments are performed for the two choices of viscosity constant, i.e., $\nu=0.1$ and $\nu=0.001$. For large viscosity constant $\nu=0.1$, submatrix $W$ is not far from a symmetric matrix since $\|S\|_2/\|H\|_2\approx 0.1272$ for $l=16$ and $0.0703$ for $l=32$. Hence, $P=\omega H$ should be a better choice than $P=(1/{\omega})(I+\omega L_s)(I+\omega U_s)$ for constraint preconditioner $M$. For small viscosity constant $\nu=0.001$, simple calculation gives $\|S\|_2/\|H\|_2\approx 12.7235$ for $l=16$ and $7.0337$ for $l=32$, which means $W$ is an essentially nonsymmetric matrix. Constraint preconditioner $M$ of the form \eqref{02} with $P=(1/{\omega})(I+\omega L_s)(I+\omega U_s)$ should be a better choice no matter for GCP iteration method \eqref{03} or for preconditioned GMRES(10) and QMR methods.

\begin{table}[!h]
\caption{Numerical results for iteration schemes \eqref{32} and \eqref{03}} \label{tab:02}
\begin{center}
\begin{tabular}{l*{9}{c}}
\hline\noalign{\smallskip}
&&&\multicolumn{3}{l}{$l=16$}&&\multicolumn{3}{l}{$l=32$}\\
\cline{4-6}\cline{8-10}\noalign{\smallskip}
&Case no.&&$\omega$&IT&CPU&&$\omega$&IT&CPU\\
\hline\noalign{\smallskip}
$\nu=0.1$& I  &&1.00 & 11 & 0.0312 && 1.00 & 8 & 0.8580 \\
& II  &&0.98 & 89 & 0.0468 && 0.99 & 249 & 1.3884 \\
& III &&$-$ & $-$ & $-$ && $-$ & $-$ & $-$ \\
& IV  &&$-$ & $-$ & $-$ && $-$ & $-$ & $-$ \\
& V &&$-$ & $-$ & $-$ && $-$ & $-$ & $-$ \\
& VI  &&$-$ & $-$ & $-$ && $-$ & $-$ & $-$ \\
$\nu=0.001$& I  &&$-$ & $-$ & $-$ && $-$ & $-$ & $-$ \\
& II  && 0.08 & 202 & 0.0624 && 0.16 & 298 & 1.5132 \\
& III  &&$-$ & $-$ & $-$ && $-$ & $-$ & $-$ \\
& IV  &&$-$ & $-$ & $-$ && $-$ & $-$ & $-$ \\
& V &&$-$ & $-$ & $-$ && $-$ & $-$ & $-$ \\
& VI  &&$-$ & $-$ & $-$ && $-$ & $-$ & $-$ \\
\noalign{\smallskip}\hline
\end{tabular}
\end{center}
\end{table}


In Table \ref{tab:02}, besides the experimentally found optimal values of parameter $\omega$, we list the iteration steps and elapsed CPU times for iteration schemes \eqref{03} and \eqref{32} with the six choices of $M$ presented in Table \ref{tab:01}. From the numerical results, we see that for large viscosity constant $\nu=0.1$, the four methods including iteration scheme \eqref{03} with $M$ defined by Cases III and IV, and iteration scheme \eqref{32} with $M$ defined by Cases V and VI, do not achieve the stop criterion \eqref{40} within IT$_{\max}=5000$ iteration steps for any $\omega\in(0,2000]$. The GCP iteration method \eqref{03} with $M$ being defined by Case I uses less iteration steps and CPU times to compute a satisfactory solution than the method with $M$ being defined by Case II. For small viscosity constant $\nu=0.001$, matrix $W$ is essentially nonsymmetric, only the GCP iteration method \eqref{03} with $M$ being defined by Case II can obtain a satisfactory solution. These numerical results are consistent with our conjecture made in the last paragraph, and also verify the robustness of the GCP iteration method.


\begin{table}[!h]
\caption{Numerical results for preconditioned GMRES(10) methods}\label{tab:03}
\begin{center}
\begin{tabular}{l*{9}{c}}
\hline\noalign{\smallskip}
&&&\multicolumn{3}{l}{$l=16$}&&\multicolumn{3}{l}{$l=32$}\\
\cline{4-6}\cline{8-10}\noalign{\smallskip}
&Case no.&&$\omega$&IT&CPU&&$\omega$&IT&CPU\\
\hline\noalign{\smallskip}
$\nu=0.1$
& I  && 1.50 &  14 & 0.0156  && 1.61  & 22 & 0.2340 \\
& II  && 0.63 & 34 & 0.0468 && 0.64  &63 & 0.5148 \\
& III  && 0.03 & 29  & 0.0312 && 0.02  & 30  & 0.4524  \\
& IV  && 0.02 & 49 & 0.0624 && 0.02  & 110  & 0.5772  \\
& V   && 0.01 & 87 &  0.0468 && 0.02  & 695  & 0.7956  \\
& VI  && $-$ & $-$  & $-$  && $-$  & $-$  & $-$ \\
$\nu=0.001$
& I  && 26.40 &  748 & 0.4056  &&  28.62 &  1340 & 1.6536  \\
& II  && 0.04 & 118 & 0.0780 && 0.05 & 700 & 0.7644 \\
& III  && 0.04 & 1846 &  0.9984 && 0.02  & 4647  & 4.5552  \\
& IV  && 0.06 &  229 & 0.2028  &&  0.10 & 851  & 1.9968  \\
& V   && 0.02 &  953 & 0.5148  &&  0.01 &  1851 & 1.9188  \\
& VI  && $-$ & $-$  & $-$  && $-$  & $-$  & $-$ \\
\noalign{\smallskip}\hline
\end{tabular}
\end{center}
\end{table}

\begin{table}[!h]
\caption{Numerical results for preconditioned QMR methods}\label{tab:04}
\begin{center}
\begin{tabular}{l*{9}{c}}
\hline\noalign{\smallskip}
&&&\multicolumn{3}{l}{$l=16$}&&\multicolumn{3}{l}{$l=32$}\\
\cline{4-6}\cline{8-10}\noalign{\smallskip}
&Case no.&&$\omega$&IT&CPU&&$\omega$&IT&CPU\\
\hline\noalign{\smallskip}
$\nu=0.1$& I  && 1.52 &  11 & 0.0624  && 1.59  & 13  & 0.3432  \\
& II  && 0.60 & 35  &  0.0936 && 0.63  & 68  &  0.5928 \\
& III  && 2.12 &  31 &  0.0936 && 2.11  & 36 &  0.6084 \\
& IV  && 1.00 & 89  & 0.0780  && 0.99  & 183  & 0.5460  \\
& V  && 1.26 &  47 & 0.0624  && 1.11  &  69 &  0.5340 \\
& VI  && 0.90 &  385 &  0.1404 && 0.85  & 1420  & 1.0296  \\
$\nu=0.001$& I  && 24.10 &  276 &  0.6084 && 21.60  & 486  &  10.0777 \\
& II  && 0.06 & 142 & 0.1560  && 0.05  &  294 & 0.6864  \\
& III  && $-$ & $-$  & $-$  &&  $-$ & $-$  & $-$  \\
& IV  && 0.09 &  217 & 0.6240  &&  0.11 &  340 &  1.9344 \\
& V  && 28.35 & 652  & 0.4836  && 25.67  & 1411  &  9.0169 \\
& VI  && 0.02  & 871  & 0.3120  &&  0.04 &  3852 &  3.2448 \\
\noalign{\smallskip}\hline
\end{tabular}
\end{center}
\end{table}

To solve the singular saddle-point problems \eqref{01}, we also use each of the choices of matrix $M$ in Table \ref{tab:01} as an preconditioner to accelerate GMRES(10) and QMR methods, respectively. The experimentally found optimal values of parameter $\omega$, iteration steps, elapsed CPU times of the preconditioned GMRES(10) and QMR methods are listed in Tables \ref{tab:03} and \ref{tab:04}. Numerical results in the two tables show that the preconditioning effects of singular constraint preconditioner $M$ defined in \eqref{02} with $P=\omega H$ (i.e., the Case I in Table \ref{tab:01}), no matter for GMRES(10) method or for QMR method, are the best for the case $\nu=0.1$. For the case $\nu=0.001$, GMRES(10) method preconditioned by singular constraint preconditioner $M$ with $P=(1/{\omega})(I+\omega L_s)(I+\omega U_s)$ (i.e., the Case II in Table \ref{tab:01}) costs the least iteration steps and CPU times comparing with other five preconditioned GMRES(10) methods, so do the preconditioned QMR methods. Hence, the singular constraint preconditioners $M$ are efficient and robust for accelerating the convergence rates of GMRES(10) and QMR methods.

\section{Conclusion}\label{sec:5}
We present a general constraint preconditioning (GCP) iteration method, i.e., iteration scheme \eqref{03} with matrix $M$ being of the form \eqref{02} and the $(1,1)$ block of $M$ being positive definite, for solving the singular saddle-point problems \eqref{01}.
The convergence properties of the GCP iteration method are carefully studied and two different choices of the $(1,1)$ block of matrix $M$ are also discussed. Theoretical analysis shows that, under suitable conditions, the GCP iteration method is convergent for any initial guess $x^{(0)}$.

Recently, based on the Hermitian and skew-Hermitian splitting (HSS) preconditioner \cite{Bai2010171,BG200520}, Zhang et al. \cite{ZhangRenZou2011273} proposed an efficient HSS-based constraint preconditioner $M$ of the form \eqref{02}, in which $P$ is chosen to be the HSS preconditioner of $W$, to solve nonsingular saddle-point problems \eqref{01}. How about the efficiency of the HSS-based constraint preconditioner used for solving singular saddle-point problems, which may be studied in our future.

\vspace{0.5cm}
\noindent{\bf Acknowledgements.} \vspace{0.2cm}

This work is partially supported by the National Basic Research (973) Program of China under
Grant 2011CB706903, the CAPES and CNPq in Brazil, the National Natural Science Foundation of China under Grant 11271174 and the Mathematical Tianyuan Foundation of China under Grant 11026064.

\bibliographystyle{model1-num-names}
\bibliography{References}
\end{document}